\documentclass[12pt]{article}
\usepackage{amssymb,latexsym}
\def\C{\mathbb C}

\def\N{\mathbb N}

\newtheorem{thm}{Theorem}[section]
\newtheorem{lem}{Lemma}[section]

\begin{document}
\sffamily
\title{Zeros of differences of meromorphic functions}
\author{Walter Bergweiler and J.K. Langley
\thanks{Research partly carried out during a visit by the second
author to the
Christian-Albrechts-Universit\"at zu Kiel in April-May 2005, 
supported by a grant from
the Alexander von Humboldt Stiftung. The first author is supported 
by the G.I.F.,
the German--Israeli Foundation for Scientific Research and
Development, Grant G -809-234.6/2003.}
}
\maketitle
\begin{abstract}
Let $f$ be a function transcendental and meromorphic in the plane, and define 
$g(z)$ by
$g(z) = \Delta f(z) = f(z+1) - f(z)$. A number of results are proved 
concerning the existence of zeros of
$g(z)$ or $g(z)/f(z)$, in terms of the growth and the poles of $f$. The results
may be viewed as discrete analogues of existing theorems on the zeros of
$f'$ and $f'/f$.
\\
{\bf MSC 2000:} 30D35.
\end{abstract} 
\section{Introduction}
Let $f$ be a function transcendental and
meromorphic in the plane. The forward differences
$\Delta^n f$ are defined in the standard way \cite[p.52]{Whi} by
\begin{equation}
\Delta f(z) = f(z+1) - f(z), \quad 
\Delta^{n+1} f(z) = \Delta^{n} f(z+1) - 
\Delta^{n} f(z) , \quad n= 0, 1, 2, \ldots .
\label{w1}
\end{equation}
This paper is concerned with the question of whether the forward
differences defined in (\ref{w1}) must have zeros, the principal
motivations for this being twofold.

First, considerable recent attention has been given to meromorphic solutions
$y = f(z)$ in the plane of difference equations
$$
a_n(z) y(z + n) + \ldots + a_1(z) y(z+1) + a_0 (z) y(z) = A(z) ,
$$
as well as of functional equations of related type. A number of 
papers (including \cite{AHH,chiangfeng,HK,HK2,HKLRT,IshYan})
focus on the growth and zeros of solutions of such equations,
investigating analogies and contrasts with the theory of linear
differential equations in the complex plane.
The second motivation is as a discrete analogue of the following theorem,
in which the notation is that of \cite{Hay2}.

\begin{thm}[\cite{ELR,Hinchliffe}]\label{elrthm}
Let $f$ be transcendental and meromorphic in the plane with
\begin{equation}
\liminf_{r\to\infty} \frac{T(r,f)}{r} = 0.
\label{elr1}
\end{equation}
Then $f'$ has infinitely many zeros. 
\end{thm}

Theorem \ref{elrthm} is sharp, as shown by $e^z, \, \tan z$ and
examples of arbitrary order greater than $1$ constructed in \cite{CER}.
The result was originally proved in \cite{ELR} 
(see also \cite{BE}) with $\limsup$ in (\ref{elr1}),
the improvement to $\liminf$ being due to Hinchliffe 
\cite{Hinchliffe}. For $f$
as in the hypotheses of Theorem \ref{elrthm} it follows from
Hurwitz' theorem that if $z_1$ is a zero of $f'$ then
$f(z+c) - f(z)$ has a zero near $z_1$, for all sufficiently small
$c \in \C \setminus \{ 0 \}$. This makes it natural to ask 
whether $f(z+c) - f(z)$, for such functions $f$, must always have
infinitely many zeros. Here there is no loss of generality in assuming that
$c=1$, since otherwise $f(z)$ may be replaced by $F(z) = f(cz)$. Examples
such as
$$
f(z) = ze^{ 2 \pi i z }  + h(z), \quad \Delta f(z) = e^{ 2 \pi i z },
$$
where $h$ is an entire function of period $1$, show that attention
must be restricted to functions of subexponential growth.

Consider first the case where $f$ is a transcendental entire function
of order less than $1$. Then so is the first difference 
$\Delta f$ \cite{Whi} (see also Lemma \ref{notrational}) and by
repetition of this argument each difference $\Delta^n f$, for $ n \geq 1$,
is transcendental entire of order less than $1$ and so obviously has 
infinitely many zeros. Thus for entire $f$ it is natural to consider
zeros not of $\Delta^n f$ but rather of the divided difference
$\frac{\Delta^n f}f$. This is analogous to the counterpart of Theorem
\ref{elrthm} for the logarithmic derivative
$f'/f$ proved in \cite{CER,Hinchliffe}: if $f$ is
transcendental entire satisfying (\ref{elr1}) or 
transcendental meromorphic with $\liminf_{r\to\infty} r^{-1/2} T(r, f)
=0$ then $f'/f$ has infinitely many zeros.

The following result may be proved using a version
of the Wiman-Valiron theory for differences developed in \cite{IshYan}: a
proof based instead on the standard Wiman-Valiron theory 
\cite{Hay5,Valiron} will be given in
\S\ref{pflesshalf}.

\begin{thm}\label{lesshalf}
Let $n \in \N$ and let
$f$ be a transcendental entire function of order $\rho < \frac12$, and set
\begin{equation}
G(z) = \frac{\Delta^n f (z)}{f(z)}  .
\label{w2}
\end{equation}
If $G$ is transcendental then $G$ has infinitely many zeros.
In particular if $f$ has order less than $\min \left\{ 
\frac1n , \frac12 \right\}$ then $G$ is transcendental
and has infinitely many zeros.
\end{thm}

The proof of Theorem \ref{lesshalf} relies upon the classical
$\cos \pi \rho$ minimum modulus theorem
\cite[Theorem 6.13, p.331]{Hay7}, and 
breaks down if the order is at least $\frac12$.  However,
for the first divided difference it is possible to extend Theorem
\ref{lesshalf} slightly beyond $\rho = \frac12$.

\begin{thm}\label{thm1}
There exists $\delta_0 \in \left(0, \frac12\right)$ with 
the following property. Let $f$ be a transcendental entire function with order
\begin{equation}
\rho (f) \leq \rho < \frac12 + \delta_0 < 1.
\label{a1}
\end{equation}
Then
\begin{equation}
G(z) = \frac{\Delta   f(z)}{f(z)} 
= \frac{f(z+1) - f(z)}{f(z)} 
\label{a2}
\end{equation}
has infinitely many zeros.
\end{thm}
It will be seen from the proof of Theorem \ref{thm1} in \S\ref{pfthm1}
that the constant $\delta_0$ is extremely small. It seems
reasonable to conjecture that the conclusion of
Theorem \ref{thm1} in fact holds for
$\rho (f) < 1$, but the present proof, which 
is based on an estimate of Miles and Rossi \cite{MR} for the size of the 
set where $f'/f$ is large, will not give this. 

In  considering the existence of zeros of 
$g(z) = f(z+c) - f(z)$ when $f$ is meromorphic, complications appear to 
arise from the poles of $f$, which may or may not be poles of $g$. 
The following theorem will be deduced in \S\ref{pfthm2} from
Theorem \ref{elrthm}, using an approximation of $g(z)$ in terms of
$f'(z)$ which will be developed in \S\ref{cartan}.

\begin{thm}\label{thm2}
Let $f$ be a function
transcendental and meromorphic of lower
order $\lambda (f) < \lambda < 1$ in the plane.
Let $c \in \C \setminus \{ 0 \}$ be such that at most finitely many poles 
$z_j , z_k $ of $f$ satisfy $z_j - z_k = c$. Then 
$g(z) = f(z+c) - f(z)$ has infinitely many zeros.
\end{thm}
It is clear that for a given $f$
all but countably many $c \in \C$ satisfy the hypotheses
of Theorem \ref{thm2},
but the following construction shows that Theorem \ref{thm2} fails without
the hypothesis on $c$, even for lower order $0$.
\begin{thm}\label{onezero}
Let $\phi (r)$ be a positive non-decreasing
function defined on $[1, \infty)$ 
which satisfies
$\lim_{r \to \infty} \phi (r) = \infty $. Then there exists
a function $f$ transcendental and meromorphic in the plane with
\begin{equation}
\limsup_{r\to\infty}\frac{T(r,f)}{r} < \infty 
\label{oz1}
\end{equation}
and
\begin{equation}
\liminf_{r\to\infty}\frac{T(r,f)}{\phi(r) \log r} < \infty , 
\label{oz2}
\end{equation}
such that
\begin{equation}
g(z) = \Delta f (z) = f(z+1) - f(z)
\label{oz3}
\end{equation}
has only one zero. Moreover, the function $g$ satisfies
\begin{equation}
\limsup_{r\to\infty}\frac{T(r,g)}{\phi(r) \log r} < \infty . 
\label{oz4}
\end{equation}
\end{thm}

On the other hand, for transcendental meromorphic functions of
sufficiently small growth, it is possible to show that either the first
difference or the first divided difference has infinitely many zeros.

\begin{thm}\label{thm3}
Let $f$ be a function transcendental and meromorphic in the plane with
\begin{equation}
T(r, f) = O( \log r )^2 \quad \hbox{as $r \to \infty$,}
\label{d1}
\end{equation}
and set 
\begin{equation}
g(z) = \Delta f(z) = f(z+1) - f(z),
\quad G(z) = \frac{\Delta f(z)}{f(z)}  = \frac{f(z+1) - f(z)}{f(z)}.
\label{d2}
\end{equation}
Then at least one of $g$ and $G$ has infinitely many
zeros.
\end{thm}

The proofs of Theorems
\ref{thm1}, \ref{thm2}, \ref{onezero} and \ref{thm3} will be given in
\S\ref{pfthm1}, \S\ref{pfthm2}, \S\ref{pfonezero} and \S\ref{pfthm3}
respectively.

\section{Preliminaries}
\begin{lem}\label{notrational}
Let $f$ be a function transcendental and meromorphic in the plane which 
satisfies $(\ref{elr1})$,
and with the notation $(\ref{w1})$ let
$g = \Delta f$ and $G = g/f$. 
Then $g$ and $G$ are both  transcendental.
\end{lem}
The assertion concerning $g$ may be found in
\cite[p.101]{Whi}, but a proof will be given for completeness.\\
{\em Proof.}
Suppose first that $G$ is a rational function. Then (\ref{a2}) gives 
\begin{equation} 
f(z+1) = R_0(z) f(z), \quad f(z-1) = R_1(z) f(z), 
\label{a3} 
\end{equation} 
where $R_0$ and $R_1$ are rational functions, neither identically zero. 
Take $r_0 > 0$, so large that $R_0$ and $R_1$ have no zeros or
poles in $|z| > r_0$. Suppose that
$z_0$ is a zero of $f$ with $|z_0| > r_0$. 
Then (\ref{a3}) shows that either $z_0 + 1, z_0 +2, \ldots $,
or $z_0 - 1, z_0 -2, \ldots $, are zeros of $f$, 
depending on the sign of ${\rm Re} \, z_0$, 
and both contradict (\ref{elr1}).
The same argument shows that $f$ cannot have a pole $z_0$ with
$|z_0| > r_0$. But (\ref{elr1}) shows that $f$ must have
infinitely many zeros or infinitely many poles,
and this is a contradiction. 

The proof that $g$ is transcendental
is similar. Assume that $g$ is a rational function. Then there
exist rational functions $R_2$ and $R_3$ such that
\begin{equation}
f(z+ 1) = f(z) + R_2(z), \quad f(z-1) = f(z) + R_3(z).
\label{a3a}
\end{equation}
If $f$ has infinitely many poles then a contradiction arises exactly as
in the proof that $G$ is transcendental. Assume henceforth that
$f$ has finitely many poles. Then there exists a rational function
$R_4$ such that $h = f - R_4$ is transcendental
entire, and by considering $h$ in
place of $f$ it may be assumed that $R_2$ is a polynomial in
(\ref{a3a}). But then by
\cite[p.21]{Whi} there exists a polynomial $P$ such that
$P(z+1) - P(z) = R_2(z)$, and so by considering $f - P$ in place
of $f$ it may now be assumed that $R_2 \equiv 0$ in
(\ref{a3a}). Hence $f$ has period $1$, which contradicts (\ref{elr1}).
\hfill$\Box$
\vspace{.1in}

A key role in the proof of Theorem \ref{thm1} will be
played by the following result of Miles and Rossi \cite{MR}.

\begin{lem}[\cite{MR}]\label{lem6}
Let $f$ be a transcendental entire function of 
order $\rho (f) \leq \rho < \infty$. Let $\gamma > 0$, and for
$r > 0$ let
\begin{equation}
U_r = \left\{ \theta \in [0, 2 \pi ] : 
\left| \frac
{ re^{i \theta } 
f'(  re^{i \theta } ) }
{f(  re^{i \theta } ) } \right| > \gamma n(r, 1/f) \right\} .
\label{mr1}
\end{equation}
Let $M > 3$. Then there exists a set $E_M \subseteq [1, \infty )$
satisfying
\begin{equation}
\underline{\rm logdens} \, E_M =
\liminf_{r\to\infty} \left( \frac1 {\log r}
\int_{[1, r] \cap E_M} \frac{dt}t \right)
\geq 1 - \frac3M ,
\label{mr2}
\end{equation}
such that 
\begin{equation}
m(U_r) > \left( \frac{1 - \gamma }{ 7M (\rho + 1) } \right)^2 
\quad \hbox{for} \quad r \in E_M ,
\label{mr3}
\end{equation}
in which $m(U_r)$ denotes the
Lebesgue measure of $U_r$ .
\end{lem}

The proof of Theorem \ref{thm1} will also
require the following variant of a standard estimate for harmonic measure
\cite[p.116-7]{tsuji}.

\begin{lem}\label{lem7}
Let $H$ be a transcendental entire function of order $\rho < \infty$.
For large $r > 0$
define $r \theta (r)$ to be the length of the longest arc of the circle
$|z| = r$
on which $|H(z)| > 1$, with $\theta(r) = 2 \pi $ 
if the minimum modulus 
\begin{equation}
m_0(r, H) = \min \{ |H(z)| : |z| = r \} 
\label{hm1}
\end{equation}
satisfies $m_0(r, H) > 1$. Then at least one of the following is
true:\\
(i) there exists a set $F \subseteq [1, \infty )$ 
of positive upper logarithmic density such that 
$m_0(r, H) > 1$ for $r \in F$;\\
(ii) for each $\tau \in (0, 1)$ the set
\begin{equation}
F_\tau = \{ r : \theta (r) > 2 \pi (1 - \tau ) \} 
\label{hm2}
\end{equation}
satisfies
\begin{equation}
\underline{\rm logdens} \, F_\tau \geq 
\frac{1 - 2 \rho ( 1 - \tau ) }\tau .
\label{hm3}
\end{equation}
\end{lem}
Note that when $\rho = \frac12$ the right hand side of (\ref{hm3}) 
is $1$, and that when $H$ has lower order less than $\frac12$ it follows
from Barry's lower order version of 
the classical $\cos \pi \rho $ theorem
\cite{BarK} (see also \cite[p.331]{Hay7}) that conclusion (i) always holds.
\\\\
{\em Proof.} Assume that conclusion (i) does not hold. Define $\theta^*(r)$
to be the same as $\theta (r)$, except that $\theta^*(r) = \infty$
if $m_0(r, H) > 1$. Then $\theta (r) = \theta^*(r)$ on a set 
of logarithmic density $1$. Since
$$
\frac1{\theta(r)} \leq 
\frac1{\theta^*(r)} +
\frac1{ 2 \pi } 
$$
for all large $r$,
the standard Carleman-Tsuji estimate for harmonic measure 
\cite[pp.116-7]{tsuji} gives
a large positive $R$ such that
\begin{equation}
\int_R^r \frac{\pi dt}{t \theta (t) } \leq 
\int_R^r \frac{\pi dt}{t \theta^* (t) } + o( \log r) \leq 
(\rho + o(1) ) \log r 
\label{hm4}
\end{equation}
as $r \to \infty$. Hence if $\tau \in (0, 1)$ then
(\ref{hm4}) leads to, as $r \to \infty$,
\begin{eqnarray*}
( 2 \rho + o(1) ) \log r 
&\geq & \int_{[R, r] \cap F_\tau } \frac{2 \pi dt}{t  \theta (t) } +
\int_{[R, r] \setminus  F_\tau } \frac{2 \pi dt}{t  \theta (t) } \\
&\geq & \int_{[R, r] \cap F_\tau } \frac{dt}t 
+ \frac1{1 - \tau } \int_{[R, r] \setminus F_\tau } \frac{dt}t \\
&=&
\int_{[R, r] \cap F_\tau } \frac{dt}t 
+  \frac1{1 - \tau } 
\left( \log r - O(1) -  \int_{[R, r] \cap F_\tau } \frac{dt}t \right) \\
&=& 
- \frac\tau{1 - \tau } \int_{[R, r] \cap F_\tau } \frac{dt}t + 
\frac1{1 - \tau }  \log r - O(1) ,
\end{eqnarray*}
from which (\ref{hm3}) follows.
\hfill$\Box$

\section{An estimate of Cartan type}\label{cartan}

Following Hayman \cite{Hay8}, define an $\varepsilon$-set to be a countable
union of discs 
\begin{equation}
E = \bigcup_{j=1}^\infty B(b_j, r_j) \quad \hbox{ such that} \quad
\lim_{j\to \infty} |b_j| = \infty \quad \hbox{and} \quad
\sum_{j=1}^\infty \frac{r_j}{|b_j|} < \infty .
\label{car1}
\end{equation}
Here and henceforth $B(a, r)$ denotes the open disc of centre $a$
and radius $r$, and $S(a, r)$ will denote the corresponding
boundary circle.
Note that if $E$ is an $\varepsilon$-set then the set of $r \geq 1$ 
for which the circle $S(0, r)$ meets $E$ has finite logarithmic 
measure and hence zero logarithmic density.

The term $\varepsilon$-set was introduced in the context
of the following theorem, which was proved by Hayman for
entire functions \cite{Hay8} and by Anderson and Clunie 
\cite{ACslow} for meromorphic functions with deficient poles.

\begin{thm}[\cite{ACslow}]\label{haythm}
Let $h$ be a function transcendental and meromorphic in the plane, with
$$
T(r, h) = O( \log r)^2  \quad \hbox{as} \quad r \to \infty ,
$$
and assume that the Nevanlinna deficiency $\delta ( \infty, h)$
of the poles of $h$ is positive.
Then there exists an $\varepsilon$-set $E$ such that
$$
\log |h(z)| \geq ( \delta ( \infty , h) - o(1) ) T(|z|, h) \quad
\hbox{as $z \to \infty$ in $\C \setminus E$}.
$$
\end{thm}

\begin{lem}\label{lem1}
Let $a_1, a_2, \ldots $ be complex numbers with
$|a_k| \leq |a_{k+1}| $ and $\lim_{k\to\infty} |a_k| = \infty$.
For $r > 0$ let $n(r)$ be the number of $a_k$, 
taking account of repetition, with $|a_k| \leq r$. Let $\alpha > 1$. 
Then there exist a positive constant
$d_\alpha$, depending only on $\alpha$, and
an $\varepsilon$-set $E = E_\alpha$ such that for large $z$ with
$z \not \in E$ and $|z| = r$,
\begin{equation}
\sum_{|a_k| < \alpha r } \frac1{|z - a_k|} <
d_\alpha \frac{n( \alpha^2 r) }r ( \log r)^\alpha \log n( \alpha^2 r) .
\label{car2}
\end{equation}
Moreover, if 
\begin{equation}
\sum_{a_k \neq 0} \frac1{|a_k|} < \infty ,
\label{car3}
\end{equation}
then for any positive constant $h$ it is possible to choose $E$ so that
$|z - a_k| \geq 2h$ for all large $z \not \in E$ and for all $k$.
\end{lem}
{\em Proof.} The first part is proved, though not explicitly stated,
by Gundersen \cite[Lemma 2]{Gun2}. In particular \cite[(5.8)]{Gun2}
shows that (\ref{car2}) holds outside an exceptional set 
satisfying (\ref{car1}).
Suppose now that (\ref{car3}) holds. Then if $k_0$ is large the set
$E' = \bigcup_{k \geq k_0} B(a_k, 2h) $ is an $\varepsilon$-set,
and it is only necessary to replace $E$ by $E \cup E'$.
\hfill$\Box$
\vspace{.1in}

The next lemma is standard and can be found, for example,
in \cite[p.65]{Jank}.

\begin{lem}[\cite{Jank}]\label{lem2}
Let $g$ be non-constant and meromorphic in the plane and let 
$\beta > 1$. Then for $|z| = r$ sufficiently large,
\begin{equation}
\left| \frac{g'(z)}{g(z)} \right| \leq d_\beta \frac{ T( \beta r , g)}r + 
\sum_{|a_k| < \beta r } \frac2{|z - a_k| } ,
\label{car4}
\end{equation}
in which $d_\beta$
is a positive constant depending only on $\beta$, and the $a_k$ are the
zeros and poles of $g$, repeated according to multiplicity.
\end{lem}

\begin{lem}\label{lem3}
Let $g$ be a function transcendental and meromorphic in the plane 
of order less than~$1$. Let $h > 0$. Then there exists an
$\varepsilon$-set $E$ such that
\begin{equation}
\frac{g'(z+c)}{g(z+c)} \to 0 \quad  \hbox{and} \quad
\frac{g(z+c)}{g(z)} \to 1 \quad
\hbox{as $z \to \infty$ in $\C \setminus E$} ,
\label{car5}
\end{equation}
uniformly in $c$ for $|c| \leq h$. Further, $E$ may be chosen so that
for large $z$ not in $E$ the function
$g$ has no zeros or poles in $|\zeta - z | \leq h$.
\end{lem}
{\em Proof.} Since $g$ has order less than $1$ the sequence $(a_k)$ of
zeros and poles of $g$, with repetition according to multiplicity, 
evidently satisfies (\ref{car3}). Apply Lemmas 
\ref{lem1} and \ref{lem2} with $\alpha = 4$ and $\beta = 2$. 
In particular, Lemma \ref{lem1} gives an $\varepsilon$-set $E$ such that
for large $z \not \in E$ the estimate     
(\ref{car2}) holds,
as well as $|z - a_k| \geq 2h$ for all $k$. 
Let $z$ be large, 
not in $E$, set $r = |z|$, and let $|\zeta - z | \leq h$. Then
\begin{equation}
| \zeta | \leq 2 r 
\quad \hbox{and} \quad 
| \zeta - a_k | \geq | z - a_k | - h 
\geq \frac{ | z - a_k | }2  \quad \forall \, k .
\label{car6}
\end{equation}
In particular, $\zeta$ is not a pole or zero of $g$.
Now (\ref{car2}), (\ref{car4})
and (\ref{car6}) give an absolute constant $d > 0$ such that
\begin{eqnarray*}
\left| \frac{g'(\zeta)}
{g(\zeta)} \right| 
&\leq& d \frac{T(4 r, g)}r + \sum_{|a_k| < 4r } \frac2{|\zeta - a_k|} \\
&\leq& d \frac{T(4 r, g)}r + \sum_{|a_k| < 4r } \frac4{|z - a_k|} \\
&\leq& d \frac{T(4 r, g)}r + 4 d_4 
\frac{n(16r)}{r} (\log r)^4 \log n(16r)
= o(1) ,
\end{eqnarray*}
where $n(r) = n(r, g) + n(r, 1/g)$.
The first assertion of (\ref{car5}) now follows immediately
on setting $\zeta = z + c$, 
while the second assertion follows on writing
$$
\log \frac{g(z+c)}{g(z)} = \int_z^{z+c} 
\frac{g'(\zeta)}{g(\zeta)} d \zeta = o(1) .
$$
\hfill$\Box$
\vspace{.1in}

The example $g(z) = \sin z$ and the remark following (\ref{car1})
show that Lemma \ref{lem3} fails for functions of order $1$, 
since for any $r > 0$ there exists $c \in (0, \pi)$ such that 
$g(r+c)/g(r)$ is either $0$ or $\infty$.

\begin{lem}\label{lem4}
Let $f$ be a function transcendental and meromorphic in the plane
of order less than~$1$. Let $h > 0$.
Then there exists an $\varepsilon$-set $E$
such that
\begin{equation}
| f(z+c) - f(z) - c f'(z) | \leq |c|^2 \frac{| f''(z) |}2 (1 + o(1) )
\quad \hbox{as $z \to \infty$ in $\C \setminus E$} ,
\label{car7}
\end{equation}
uniformly in $c$ for $|c| \leq h$.
\end{lem}
{\em Proof.} Apply Lemma \ref{lem3}, with $g = f''$. This gives an
$\varepsilon$-set $E$ satisfying the second assertion of
(\ref{car5}) with $g = f''$, and
for large $z $ not in $E$ there are no zeros or 
poles of $f''$ in $|\zeta - z| \leq h$. 

Let $z$ be large, not in $E$. Then 
(\ref{car5}) gives, for $|u| \leq h$,
$$
| f'(z+u) - f'(z) | \leq \int_z^{z+u} | f''( \zeta ) | \, 
| d \zeta | \leq 
|u f''(z)| (1 + o(1)) 
$$
and so, for $|c| \leq h$,
$$
| f(z+c) - f(z) - c f'(z) | = 
\left| \int_0^c ( f'( z + u ) - f'(z) ) du \right| \leq 
(1 + o(1)) | f''(z)| \int_0^{|c|} t dt ,
$$
from which (\ref{car7}) follows at once.
\hfill$\Box$
\vspace{.1in}

\begin{lem}\label{lem5}
Let $f$ and $h$ be as in Lemma $\ref{lem4}$. Then there exists an
$\varepsilon$-set $E'$ such that 
\begin{equation}
f(z + c) - f(z) = c f'(z) (1 + o(1)) \quad
\hbox{as $z \to \infty$ in $\C \setminus E'$} ,
\label{car8}
\end{equation}
uniformly in $c$ for $|c| \leq h$.
\end{lem}
{\em Proof.} 
Lemma \ref{lem5} follows immediately from 
Lemma \ref{lem4}, it being only necessary to adjoin to the 
$\varepsilon$-set $E$ of Lemma \ref{lem4} an
$\varepsilon$-set $E''$ outside which
$f''(z)/f'(z) \to 0$, which is possible by Lemma \ref{lem3}.
\hfill$\Box$
\vspace{.1in}

For functions of lower order less than $1$ the condition (\ref{car3})
may fail, but the following weaker
assertion will suffice for subsequent application. 

\begin{lem}\label{lemlo}
Let $f$ be a  function transcendental and meromorphic in the 
plane of lower order $\lambda (f) < \lambda < 1$.
Then there exist arbitrarily large $R$ with the following
properties. First,
\begin{equation}
T(32R, f') < R^\lambda  .
\label{c2}
\end{equation}
Second, there exists
a set $J_R \subseteq [R/2, R    ]$ of linear measure 
$(1 - o(1) ) R/2$ 
such that, for $r \in J_R$,
\begin{equation}
f(z+1) - f(z) \sim f'(z) \quad \hbox{on}  \quad |z| = r .
\label{c1}
\end{equation}
\end{lem}
{\em Proof.} Let $(a_k)$ be the sequence of all 
zeros and poles of $f'$, 
with repetition according to 
multiplicity, and let $n(r)$ be the counting function of the 
$a_k$ as in Lemma \ref{lem1}. Let
$E$ be the $\varepsilon$-set arising from Lemma \ref{lem1}, with the choice
$\alpha = 4$.

It is clear from the hypotheses that
there exist arbitrarily large $R$ satisfying 
(\ref{c2}). For such $R$
let $E_R$ be the union of discs given by 
\begin{equation}
E_R = \bigcup_{|a_k| \leq  4R} B(a_k, 2) .
\label{c3}
\end{equation}
Then by (\ref{c2}) and the remark following
(\ref{car1}) there exists a subset $J_R$ of $[R/2, R]$, 
of measure $(1 - o(1)) R/2$, such that for $r \in J_R$ the
circle $S(0, r)$ does not meet $E \cup E_R$. 

Let $|z| = r \in J_R$ and let $|\zeta - z | \leq 1$. Then 
$| \zeta | \leq 2r $ and
$| \zeta - a_k | \geq \frac12 | z - a_k | $ for $|a_k| \leq 4R$, by
(\ref{c3}).
Thus Lemma \ref{lem2} with $\beta = 2$, (\ref{car2})
and (\ref{c2}) give, for
some positive constants $c_j$, 
\begin{eqnarray}
\left| \frac{f''(\zeta)}{f'(\zeta)} \right| 
&\leq & c_1 \frac{ T(4r, f')}r + 
\sum_{|a_k| < 4r } \frac2{|\zeta - a_k | } \nonumber \\
&\leq & c_2 \frac{T(4r, f')}r + 
\sum_{|a_k| < 4r } \frac4{|z - a_k | } \nonumber \\
&\leq & c_2 \frac{T(4r, f')}r + 
c_3 \frac{n(16r)}r (\log r)^4 \log n(16r) \nonumber \\
&=& o(1) .
\label{c4}
\end{eqnarray}
For $|\zeta - z| \leq 1$ and $|z| = r \in J_R$, integration of
(\ref{c4})
now leads to
$$
f'(\zeta ) \sim f'(z) , \quad
f(z+1) - f(z) 
= \int_z^{z+1} f'(\zeta )  d \zeta 
= \int_z^{z+1} f'(z) (1 + o(1)) \, d \zeta 
\sim   f'(z) 
$$
which gives (\ref{c1}). 
\hfill$\Box$
\vspace{.1in}

\noindent
{\em Remark.} 
The papers \cite{chiangfeng} and \cite{HK} include independently
obtained estimates for
the proximity function $m(r, g(z+c)/g(z) )$, when $g$ is a
meromorphic function of finite order. Applications of these estimates
appear in \cite{chiangfeng,HK,HK2}. The paper \cite{chiangfeng}, of which
the authors became aware after writing this paper, also
contains pointwise estimates for the modulus $|g(z+c)/g(z)|$ outside an 
$\varepsilon$-set, obtained via the Poisson-Jensen formula
and valid for meromorphic $g$ of arbitrary growth. However
for the applications of the present paper it is necessary to show as in
(\ref{car5}) that the function
$g(z+c)/g(z)$ itself, rather than just its modulus,
tends to $1$ outside an $\varepsilon$-set.

\section{Higher differences}\label{higher}
The aim of this section is to prove an asymptotic formula for the
higher differences $\Delta^n f$, for $ n \geq 2$, when $f$ is a transcendental
meromorphic function in the plane of order less than $1$. It will be
convenient to write
\begin{equation}
g_n(z) = \Delta^n f(z), \quad n \in \N.
\label{hd1}
\end{equation}
\begin{lem}\label{gnlem}
With the notation $(\ref{w1})$ and
$(\ref{hd1})$,
\begin{equation}
g_n'(z) = (\Delta^n f')(z) , \quad n \in \N.
\label{hd2}
\end{equation}
\end{lem}
{\em Proof.} The relation (\ref{hd2}) for $n=1$ follows immediately
on writing
$$
(\Delta f')(z) = f'(z+1) - f'(z) = g_1'(z).
$$
Assume now that $m \in \N$ and that (\ref{hd2}) is true for $1 \leq n \leq m$.
Then (\ref{w1}) gives
\begin{eqnarray*}
(\Delta^{m+1} f') (z) &=&
(\Delta^m f') (z+1) - 
(\Delta^m f') (z) \\
&=& g_m'(z+1) - g_m'(z) \\
&=&
(\Delta g_m')(z) = 
(\Delta g_m)'(z) = (\Delta^{m+1} f)'(z).
\end{eqnarray*}
\hfill$\Box$
\vspace{.1in}

\begin{lem}\label{hdlem}
Let $n \in \N$.
Let $f$ be transcendental and meromorphic of order less than $1$ in the
plane. Then there exists an $\varepsilon$-set $E_n$ such that
\begin{equation}
\Delta^n f (z) \sim f^{(n)}(z) \quad \hbox{as $z \to \infty $ in 
$\C \setminus E_n$}.
\label{hd3}
\end{equation}
\end{lem}
{\em Proof.} For $n =1$ the conclusion 
(\ref{hd3}) follows at once from (\ref{w1}) and Lemma \ref{lem5}.
Assume now that $n \in \N$ and that the lemma has been proved for $n$.
Then $g_n$ is a transcendental meromorphic function of order 
less than $1$, by Lemma \ref{notrational}, and so
there exists an $\varepsilon$-set $F_0$ such that
$$
\Delta^{n+1} f(z) = (\Delta g_n)(z) \sim g_n'(z) 
\quad \hbox{as $z \to \infty$ in $\C \setminus F_0$} .
$$
Since $f'$ also has order less than $1$ the induction hypothesis
gives an $\varepsilon$-set $F_n$ such that
$$
(\Delta^{n} f')(z)  \sim f^{(n+1)}(z) 
\quad \hbox{as $z \to \infty$ in $\C \setminus F_n$} .
$$
But $E_{n+1} =
F_0 \cup F_n$ is an $\varepsilon$-set and so the result for $n+1$
follows using (\ref{hd2}).
\hfill$\Box$
\vspace{.1in}

\section{Proof of Theorem \ref{lesshalf}}\label{pflesshalf}

Let $n \in \N$ and let $f$ be a transcendental entire function of 
order less than $\frac12$. Let $G$ be defined by (\ref{w2}). Then
Lemma \ref{hdlem} gives an $\varepsilon$-set $E_n$ such that (\ref{hd3}) holds.
Since $f$ is transcendental entire the Wiman-Valiron theory
\cite{Hay5,Valiron} may be applied to $f$. Let
$$
f(z) = \sum_{k=0}^\infty a_k z^k 
$$
be the Maclaurin series of $f$. For $r > 0$ the maximum term
$\mu (r)$ and central index $N(r)$ are defined by
$$
\mu (r) = \max \{ |a_k| r^k : k = 0, 1, 2, \ldots \} , \quad
N(r) = \max \{ k: |a_k| r^k  = \mu (r) \} . 
$$
The Wiman-Valiron theory \cite{Hay5,Valiron}
then gives a subset $F_0$ of $[1, \infty )$ of finite logarithmic
measure such that, for large $r$ not in $F_0$,
\begin{equation}
\frac{f^{(n)}(z)}{f(z)} 
\sim 
\frac{N(r)^n}{z^n} \quad \hbox
{for all $z$ satisfying} \quad 
|z| = r , \quad |f(z)| = M(r, f).
\label{lh2}
\end{equation}
By the remark following (\ref{car1}) it may be assumed 
that for large $r$ not in $F_0$ the circle $S(0, r)$ does not
meet the $\varepsilon$-set $E_n$ of (\ref{hd3}). Combining 
(\ref{hd3}) and (\ref{lh2}) then gives, for large $r$ not
in $F_0$,
\begin{equation}
G(z) \sim 
\frac{N(r)^n}{z^n} \quad \hbox
{for all $z$ satisfying} \quad 
|z| = r , \quad |f(z)| = M(r, f).
\label{lh1}
\end{equation}
If $f$ has order of growth $\rho < 1/n$ it follows from
(\ref{lh1}) that $G$ cannot be a rational function, since 
standard results from the Wiman-Valiron theory \cite{Hay5,Valiron}
imply that
\begin{equation}
\lim_{r \to \infty} N(r) = \infty , \quad
\limsup_{r \to \infty} \frac{\log N(r)}{\log r}  = \rho ,
\label{Nrho}
\end{equation}
so that $N(r)^n$ tends to infinity with $N(r)^n = o( r )$.

Assume henceforth that $G$ is transcendental, but has finitely many
zeros. Then $1/G$ is transcendental of order less than
$\frac12$ with finitely many poles. The classical 
$\cos \pi \rho$ theorem \cite[p.331]
{Hay7} now gives a positive constant $c_1$ and
a subset $F_1$ of $[1, \infty )$, of positive lower logarithmic density,
such that for large $r \in F_1$, 
\begin{equation}
\log |G(z)| < - c_1 T(r, G) < - n \log r \quad \hbox{on $|z| = r$},
\label{lh3}
\end{equation}
using the fact that $G$ is transcendental.
Since $F_0$ has finite logarithmic measure it may be assumed
without loss of generality that
$F_1 \cap F_0$ is empty, so that (\ref{lh1}) holds for large
$r \in F_1$. But (\ref{Nrho}) shows that
(\ref{lh1}) and (\ref{lh3}) are incompatible,
and this contradiction completes the proof of Theorem \ref{lesshalf}.
\hfill$\Box$
\vspace{.1in}

\section{Proof of Theorem \ref{thm1}}\label{pfthm1}
Let $f$ be a transcendental entire function of order 
$\rho < 1$, let $G$ be defined by (\ref{a2}),
and assume that $G$ has finitely many zeros. 
By Lemma \ref{notrational}, the function
$G$ is transcendental.

Since $\rho (f) < 1$,  Lemmas \ref{lem3} and \ref{lem5} give
a set $G_0 \subseteq [1, \infty )$ of logarithmic density $1$ 
such that 
\begin{equation} 
G(z) \sim \frac{f'(z)}{f(z)} = o(1) \quad 
\hbox{as $z \to \infty$ with $|z| \in G_0$} .  \label{a4} 
\end{equation}
Since $G$ has finitely many zeros by assumption there exists
a rational function $R_0$ with $R_0(\infty)$ finite such that
\begin{equation}
H(z) = \frac1{2z} \left( \frac1{G(z)} - R_0(z) \right) 
\label{a5}
\end{equation}
is entire and transcendental, of order at most $\rho$, and there
exists $r_1 > 0$ such that
\begin{equation}
|G(z)| <  \frac1{|z|} \quad 
\hbox{ for $|z| \geq r_1, \, |H(z)| > 1$}.
\label{a6}
\end{equation}
Apply Lemma \ref{lem7} to $H$, and suppose first that 
conclusion (i) of that
lemma holds, so that there exists a set $J_\delta \subseteq [1, \infty )$,
of positive upper logarithmic density $\delta$, on which the 
minimum modulus $m_0(r, H)$ exceeds $1$, where 
$m_0(r, H)$ is defined by (\ref{hm1}). There is no loss of 
generality in assuming that $J_\delta \subseteq G_0$, where $G_0$ 
is as in (\ref{a4}). Let $\gamma$ be small and positive,
and apply Lemma \ref{lem6} to $f$. Then (\ref{a4}) and (\ref{a6}) show
that the set §$U_r$ as defined in (\ref{mr1}) is empty for large 
$r \in J_\delta $, so that with $E_M$ as defined in Lemma \ref{lem6}
the intersection $E_M \cap  J_\delta $ is bounded, for any choice of $M > 3$,
by (\ref{mr3}). 
Since $M$ may be chosen so large that $1/3M < \delta$, this
contradicts (\ref{mr2}).

Assume henceforth that $H$ satisfies conclusion (ii) of Lemma \ref{lem7}. 
Let $M > 3$ and again let $\gamma$ be small and
positive, and define $\tau$ by
\begin{equation}
2 \pi \tau = 
\left( \frac{1 - \gamma }{7 M ( \rho + 1) } \right)^2 .
\label{a7}
\end{equation}
Let $F_\tau$ and $\theta (r)$ be as in Lemma \ref{lem7},
and again apply Lemma \ref{lem6} to $f$. This gives a subset $E_M$ of
$[1, \infty)$ satisfying (\ref{mr2}) and (\ref{mr3}), and there is
no loss of 
generality in assuming that $E_M \subseteq G_0$, where $G_0$ 
is as in (\ref{a4}).
But (\ref{mr1}), (\ref{hm2}), (\ref{a4}), (\ref{a6}),
(\ref{a7}) and the definition of $\theta (r)$
show that the intersection $E_M \cap F_\tau$ is bounded,
which by (\ref{mr2}) and (\ref{hm3}) forces
$$
1 - 2 \rho (1 - \tau ) \leq \frac{3 \tau }M   
$$
and hence
$$
2 \rho - 1 \geq \frac{\tau}{1 - \tau} \left(  1 - \frac{3 }M \right)
\geq \tau \left(  1 - \frac{3 }M \right).
$$
Since $\rho < 1$ and $\gamma$ is small, it follows using (\ref{a7})
that $\rho$ must satisfy
$$
2 \rho - 1 \geq \frac1{2 \pi } \left( \frac1{14M} \right)^2
\left( 1 - \frac3{M} \right) = h(M).
$$
In the last inequality the right hand side $h(M)$ has a maximum relative
to the interval $(3, \infty)$ at $M = 9/2$, with 
$h(9/2) = 1/23814 \pi $.
\hfill$\Box$
\vspace{.1in}

\section{Proof of Theorem \ref{thm2}}\label{pfthm2}

Let $f$ and $c$ be as in the hypotheses.
There is no loss of generality in assuming that $c = 1$. 
By the hypotheses there exist arbitrarily large $R$ satisfying the
conclusions of Lemma \ref{lemlo}.
For such $R$ let
\begin{equation}
E_R = \{ r \in  [R/2, R]  : n(r, f) = n(r-1, f) \} .
\label{b1}
\end{equation}
Then $E_R$ has linear measure 
\begin{equation}
m(E_R) \geq (1 - o(1)) R/2.
\label{b2}
\end{equation}
To see this, note that there are at most $o(R)$
points $s_k \in [R/4, R]$ at which $n(t, f)$ is discontinuous, by
(\ref{c2}). But if $r \in [R/2, R]$ is such that
$n(r) > n(r-1)$ then $r \in [s_k , s_k + 1 ]$ for some $k$.
This proves (\ref{b2}). 

Since $R$ satisfies the conclusions of
Lemma \ref{lemlo}, it follows using
(\ref{b2}) 
that there exists $r $ in $ E_R \cap J_R$ 
such that (\ref{c1}) holds, and such that 
$f(z), f(z+1)$ and $f'(z)$ have no zeros or poles on $|z| = r$. 
But by the hypotheses there exists $r_0 > 0$, independent of $R$ and $r$,
such that if $f$ has a pole of multiplicity $m$ at $z_0$ and
$r_0 \leq |z_0| \leq r-1 $ then $g(z)$ has poles at
$z_0$ and $z_0 - 1$ of multiplicity $m$.
Thus (\ref{c1}) and Rouch\'e's theorem give
\begin{eqnarray*}
n(r, 1/g) &=& n(r, 1/f') - n(r, f') + n(r, g) \\
&\geq&  n(r, 1/f') - n(r, f') + 2 n(r-1, f) - O(1)\\
&=&
n(r, 1/f') - n(r, f') + 2 n(r, f) - O(1) \\
&\geq&  n(r, 1/f') - O(1),
\end{eqnarray*}
and the result now follows since $f'$ has infinitely many zeros
by Theorem \ref{elrthm}.
\hfill$\Box$
\vspace{.1in}

\section{Proof of Theorem \ref{onezero}}\label{pfonezero}

Let $n_1, n_2, \ldots $ be positive integers with
\begin{equation}
\lim_{k\to\infty} \frac{n_{k+1}}{n_k} = \infty .
\label{oz5}
\end{equation}
Let 
\begin{equation}
H(z) = \prod_{k=1}^\infty \left( 1 + \frac{z}{A_k} \right),
\quad A_k = 4 n_k^4 .
\label{oz6}
\end{equation}
Then (\ref{oz5}) shows that
$H$ is an entire function with $T(r, H) = O( \log r )^2$ so that
by Theorem \ref{haythm} there exists an $\varepsilon$-set $E$ such that
\begin{equation}
\log |H(z)| \geq (1 - o(1))  T( |z| , H) \quad \hbox{as $z \to \infty $ with
$z \not \in E$}.
\label{oz7}
\end{equation}
Let $k$ be large and set $H_k(z) = H(z)/(z + A_k)$. 
Then (\ref{oz5}) and (\ref{oz7}) imply that there exists 
$d_k \in (2, 3)$ such that
$$
\lim_{k\to\infty} \frac{\log m_k}{\log A_k} = + \infty ,
\quad m_k = \min \{ |H_k(z)| : |z| = d_k^{\pm 1} A_k \} .
$$
In particular, since $H_k(- A_k) = H'(-A_k) $ and
$H_k$ has no zeros in $d_k^{-1} A_k \leq |z| \leq d_k A_k$ 
it follows using
(\ref{oz5}) and (\ref{oz6}) and the minimum principle that
\begin{equation}
\sum_{k=1}^\infty \left| \frac{A_k^N}{ H'(-A_k) } \right| < \infty 
\label{oz8}
\end{equation}
for any choice of $N > 0$.

Let
\begin{equation}
h(z) = \frac{H(z^4)}z .
\label{oz9}
\end{equation}
Then it follows from (\ref{oz6}) that $h$ has zeros at the points
$\pm n_k \pm i n_k$. Also if $\beta$ is a zero of $h$ then so is
$i \beta$ and
\begin{equation}
h'(\beta) = 4 \beta^2 H'(\beta^4) \neq 0, \quad
h'(i \beta) = - h'( \beta) . 
\label{oz10}
\end{equation}
By (\ref{oz6}), (\ref{oz8}), (\ref{oz9}) and (\ref{oz10}),
\begin{equation}
\sum_{k=1}^\infty n_k |c_k| < \infty , \quad
c_k = \frac1{h'(-n_k+in_k)} .
\label{oz11}
\end{equation}

Set
\begin{equation}
g(z) = \sum_{k=1}^\infty c_k \left[
\left( \frac1{z + n_k - i n_k } - \frac1{z - n_k - i n_k } \right) -
\left( \frac1{z + n_k + i n_k } - \frac1{z - n_k + i n_k } \right) 
\right] .
\label{oz12}
\end{equation} 
The series in (\ref{oz12}) converges absolutely and uniformly in 
each bounded region of the plane, by (\ref{oz11}), and the function $g$ is
meromorphic in the plane. By (\ref{oz10}) and (\ref{oz11}), 
the function $G = g - 1/h$ is entire. But (\ref{oz5}),
(\ref{oz7}) and (\ref{oz11}) imply that there exist $c > 0$ and  $n_k' \in
(4 n_k , 8 n_k )$ with
$$
M(  n_k', G) \leq o(1) + 
M(  n_k', g) \leq o(1) + c \sum_{j=1}^\infty \frac{ |c_j| }{n_k} = o(1) ,
$$
and so $G \equiv 0$ and $g$ has one zero, by (\ref{oz9}). 

Finally, set
\begin{equation}
f(z) = \sum_{k=1}^\infty c_k \left[
\sum_{j=-n_k}^{n_k-1} \frac1{z+j-in_k} 
- \sum_{j=-n_k}^{n_k-1} \frac1{z+j+in_k} \right] ,
\label{oz13}
\end{equation}
the series again convergent by (\ref{oz11}). Then $f$ and $g$ satisfy
(\ref{oz3}). A result of Keldysh \cite[p.327]{GO} (see also
\cite{CER,ELR}) gives
\begin{equation}
m(r, f) + m(r, g) = o(1) \quad \hbox{as} \quad r \to \infty .
\label{oz14}
\end{equation}
But (\ref{oz5}), (\ref{oz12}) and (\ref{oz13}) give 
\begin{equation}
n(r, f)  = O( n_k) =  O(r) \quad \hbox{and} \quad n(r, g) = O(k)
\quad \hbox{for} \quad
\sqrt{2}  n_k  \leq r < 
\sqrt{2}  n_{k+1} .     
\label{oz15}
\end{equation}
Hence (\ref{oz1}) follows using (\ref{oz14}), and since the sequence
$(n_k)$ may be chosen to grow arbitrarily fast in 
(\ref{oz5}), 
applying (\ref{oz14}) and (\ref{oz15}) again gives
(\ref{oz2})  and (\ref{oz4}).
\hfill$\Box$
\vspace{.1in}

\section{Proof of Theorem \ref{thm3}}\label{pfthm3}

Assume that $f$, $g$ and $G$ are as in the hypotheses, but that $G$ has
finitely many zeros. By Lemma
\ref{notrational}, the function $G$ is transcendental, and 
$T(r, G) = O( \log r)^2$ as $r \to \infty$,
using (\ref{d1}). By Lemma \ref{lem5} and
Theorem \ref{haythm}, there exists an $\varepsilon$-set $E$ such that
\begin{equation}
G(z) \sim \frac{f'(z)}{f(z)} \quad \hbox{and} \quad
\log |G(z)| \leq ( - 1 + o(1)) T(r, G) \quad \hbox{for $z \not \in E$ and
$|z| = r$ large.}
\label{d3}
\end{equation}
Choose $t \in [0, 2 \pi ]$ such that the ray $\arg z = t$ has bounded
intersection with $E$. Let $r_0$ be large and positive.
Integrating $f'/f$ using (\ref{d3}) along the ray 
$z = r e^{it} , r \geq r_0$,  and then around 
circles $S(0, r)$ which do not intersect $E$ then shows that there
exist a constant $b \in \C \setminus \{ 0 \}$ and 
a set $E_0 \subseteq [1, \infty )$ of finite logarithmic measure
such that
\begin{equation}
f(z) = b + o(1) \quad \hbox{for $|z| = r \in [1, \infty) 
\setminus E_0.$}
\label{d4}
\end{equation}
Set 
\begin{equation}
F(z) = f(z) - b, \quad H(z) = 
\frac{\Delta F(z)}{F(z)}
= \frac{\Delta f(z)}{f(z) - b} \, ,
\label{d5}
\end{equation}
and assume that $H$ has finitely many zeros. Then the same 
reasoning as above shows that $H$ is transcendental and that
there exists a non-zero constant $d$ such that $F(z) \sim d$ 
for $|z| = r$ large and lying outside a set of finite logarithmic
measure. This contradicts (\ref{d4}), and so $H$ must have
infinitely many zeros. 

Let $z_0$ be a zero of $H$ with $|z_0|$ large. Then $z_0$ is not
a pole of $f$, because otherwise the formula
$$
G(z) = \frac{\Delta f(z)}{f(z)} = H(z) 
\frac{f(z) - b}{f(z)} 
$$
shows that $z_0$ is a zero of $G$, which contradicts the assumption
that $G$ has finitely many zeros.
It now follows from (\ref{d5}) that
$z_0$ is a zero of $\Delta f$, and Theorem \ref{thm3} is proved.
\hfill$\Box$
\vspace{.1in}

\noindent
{\em Remark.} It seems highly unlikely that the hypothesis (\ref{d1}) in
Theorem \ref{thm3} is sharp. However the $\varepsilon$-set $E'$ arising from
Lemma \ref{lem5} may be reasonably large, at least locally, so that for
$f$ of larger growth than (\ref{d1}) difficulties may arise in integrating
$f'/f$ on the set where $G$ is small.

\noindent
Mathematisches Seminar, Christian-Albrechts-Universit\"at zu Kiel,\\
Ludewig-Meyn-Str. 4,
D-24098 Kiel, Germany.\\
bergweiler@math.uni-kiel.de\\\\
School of Mathematical Sciences, University of Nottingham, NG7 2RD, UK.\\
jkl@maths.nott.ac.uk

\end{document}